\documentclass[a4paper,12pt]{article}
\usepackage{cite,amsmath,amssymb,amsthm}
\usepackage[top=2.5cm,bottom=2.5cm,left=2.5cm,right=2.5cm]{geometry}
\usepackage[margin=1cm,%
           font=small,%
           format=hang,%
           labelsep=period,%
           labelfont=bf]{caption}
\usepackage{graphicx}

\newtheorem{proposition}{Proposition}

\newtheorem{problem}{Problem}
\newtheorem{conjecture}{Conjecture}
\begin{document}

\baselineskip=0.3in

\vspace*{3\baselineskip}

\vspace*{2mm}

\begin{center}
\textbf{\Large On the Wiener complexity and the Wiener index of fullerene graphs}
\end{center}

\begin{center}
{\large Andrey A. Dobrynin$^{1,2}$, Andrei Yu. Vesnin$^{1,2,3}$}
\end{center}

\baselineskip=0.2in

\begin{center}
\emph{$^1$Novosibirsk State University, Novosibirsk, 630090, Russia \\
\smallskip
$^2$Sobolev Institute of Mathematics, Siberian Branch of the \\
Russian Academy of Sciences, Novosibirsk, 630090, Russia\\
\smallskip
$^3$Tomsk State University, Tomsk, 634050, Russia \\
\smallskip
{\rm dobr@math.nsc.ru, vesnin@math.nsc.ru}
}
\end{center}

\begin{center}
\  %(Received May   , 2019)
\end{center}

\vspace{1mm}

\baselineskip=0.25in

\begin{abstract}
Fullerenes are molecules in the form of cage-like polyhedra,
consisting solely of carbon atoms.
Fullerene graphs are mathematical models of fullerene molecules.
The transmission of a vertex $v$ of a graph is the sum of distances from $v$ to all
the other vertices.
The number of different vertex transmissions is called the
Wiener complexity of a graph.
Some calculation results on the Wiener complexity and the Wiener index of fullerene graphs of order
$n \le 216$ are presented.
Structure of graphs with the maximal Wiener complexity or the maximal Wiener index is discussed
and formulas for the Wiener index of several families of graphs are obtained.
\end{abstract}

\section{Introduction}

A fullerene is a spherically shaped molecule consists of carbon atoms
in which every carbon ring is either a pentagon or a hexagon
and every atom has  bonds  with  exactly  three  other  atoms.
The molecule may be a hollow sphere, ellipsoid, tube, or many other shapes and sizes.
Fullerenes have been the subject of intense research, both for their chemistry and for their technological applications,
especially in nanotechnology and materials science
\cite{Ashr16,Cata11}.

Molecular graphs of fullerenes are called \emph{fullerene graphs}.
A fullerene graph is a \mbox{3-connected} 3-regular planar graph with only pentagonal
and hexagonal faces.
By Euler's polyhedral formula, the number of pentagonal faces is always 12.
It is known that fullerene graphs with $n$ vertices exist for all even $n \ge 24$ and for $n = 20$.
The number of all non-isomorphic fullerene graphs can be found in
\cite{Brin97,Fowl95,Goed15-2}.
The set of fullerene graphs with $n$ vertices will be denoted as $F_n$.
The number of faces of graphs in $F_n$ is $f = n/2 + 2$ and, therefore,
the number of hexagonal faces is $n/2-10$.
Despite the fact that
the number of pentagonal faces is negligible compared to the number of hexagonal faces,
their location is crucial to the shape and properties of fullerene molecules.
Fullerenes where no two pentagons are adjacent, \emph{i.\,e.}, each pentagon is surrounded
by five hexagons, satisfy the isolated pentagon rule and
called \emph{IPR fullerene}.
The number of all non-isomorphic IPR fullerenes was reported, for example, in
\cite{Goed15-1,Goed15-2}.
They are considered as  thermodynamic stable fullerene compounds.
Description of mathematical properties of fullerene graphs can be found in
\cite{Ando16,Ashr16,Cata11,Fowl01,Fowl95,Schw15}.

The vertex  set of a graph $G$ is denoted by $V(G)$.
The number of vertices of $G$ is called its \emph{order}.
By distance $d(u,v)$ between vertices
$u,v \in V(G)$ we mean the standard distance of a simple graph
$G$, \emph{i.\,e.}, the number of edges on a shortest path connecting
these vertices in $G$.
The maximal distances between vertices of a graph $G$ is called the \emph{diameter}
$D(G)$ of $G$.
Vertices are \emph{diametrical} if the distance between them is equal to
the diameter of a graph.
The \emph{transmission} of vertex $v \in V(G)$ is defined as the sum
of distances from $v$ to all the other vertices of $G$,
$tr(v)=\sum_{u\in V(G)} d(v,u)$.
Transmissions of vertices are used for design of many distance-based topological indices
\cite{Shar20}.
Usually, a topological index is a graph invariant that maps
a set of graphs to a set of numbers such that invariant values coincide for isomorphic
graphs.
A half of the sum of vertex transmissions gives the \emph{Wiener index}
that has found important applications in chemistry
(see selected books and reviews
\cite{Dehm14,Dobr01,Dobr02,Gutm12-1,Gutm12-2,Gutm86,Knor16,Tode00,Trin83}),
$$
W(G) = \sum_{\{u,v\}\subseteq V(G)} d(u,v) = \frac{1}{2}\sum_{v \in V(G)} tr(v).
$$

The Wiener index was introduced as structural descriptor for acyclic
organic molecules by Harold Wiener
\cite{Wien47}.
The definition of the Wiener index in terms of distances between
vertices of a graph was first given by Haruo Hosoya
\cite{Hoso71}.

The number of different vertex transmissions in a graph $G$ is known as the \emph{Wiener complexity}
\cite{Aliz16}
(or the \emph{Wiener dimension}  \cite{Aliz14}), $C_W(G)$.
This graph invariant is a measure of transmission variety.
A graph is called \emph{transmission irregular} if it has the largest possible Wiener complexity over all graphs
of a given order, \emph{i.\,e.}, vertices of the graph have pairwise different transmissions.
Various properties of transmission irregular graphs were studies in
\cite{Aliz16,Aliz18,Klav18}.
It was shown that almost all graphs are not transmission irregular.
Infinite families of transmission irregular graphs were constructed
for trees, 2-connected graphs and 3-connected cubic graphs  in
\cite{Aliz18,Dobr18,Dobr19-1,Dobr19-2,Dobr19-3}.

In this paper, we present some results of studies of
the Wiener complexity and the Wiener index of fullerene graphs.
In particular, we are interested in two
questions: does there exist of a transmission irregular fullerene graph
and can a graph with the maximal Wiener complexity has the maximal Wiener index?

\section{Distribution of graphs with respect to their Wiener complexity}

\begin{figure}[h]
\center{\includegraphics[width=\linewidth]{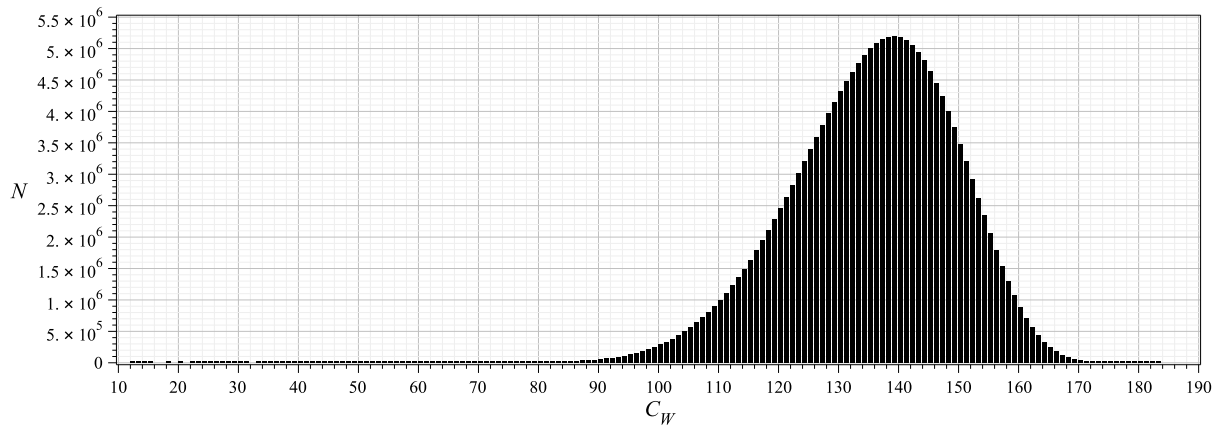}}
\caption{Distribution of fullerene graphs of $F_{196}$ with respect to
their Wiener complexity ($N$ is the number of graphs).}
\label{Fig1}
\end{figure}

Distributions of fullerene graphs with respect to their Wiener complexity have been
obtained for $n \le 216$ vertices ($f \le 110$ faces).
As an illustration, we present data for graphs with 196 vertices (100 faces).
The number of graphs of $F_{196}$ is 177\,175\,687.
Distribution of graphs of this family with respect to $C_W$ is presented in
Table~\ref{TDistr100}. A graphical representation of these data is shown in
Fig.~\ref{Fig1}.

\section{Wiener complexity of fullerene graphs}

Denote by $C_n$ the maximal Wiener complexity among all fullerene graphs
with $n$ vertices,
\emph{i.\,e.}, $C_n = \max \{C_W(G)\, |\, G \in F_n \}$.
Fullerene graphs with maximal Wiener complexity have been examined for
$n \le 216$ vertices.
Let $g_n$ be a difference between order and the Wiener complexity,
 $g_n=n-C_n$.
Then a transmission irregular graph has $g_n=0$.
It is obvious that a transmission irregular graph
has the identity automorphism group.

\begin{figure}[h!]
\center{\includegraphics[width=\linewidth]{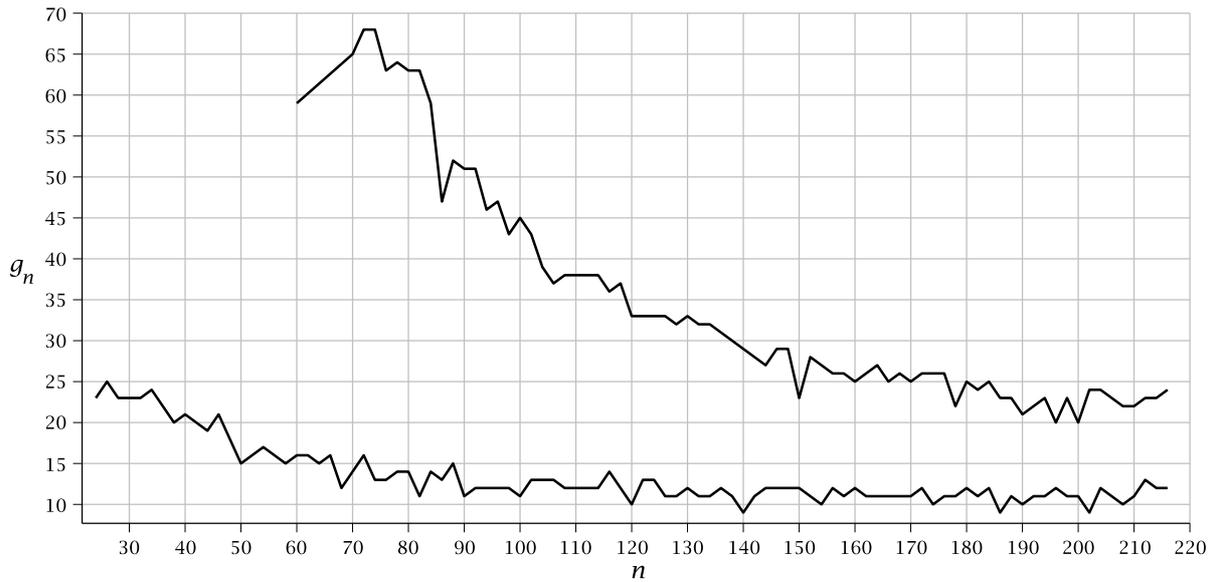}}
\caption{Difference $g_n$ between order and the maximal Wiener complexity of \\ fullerene graphs.}
\label{Fig2}
\end{figure}

The behavior of $g_n$ when the number of vertices $n$ increases
is shown in Fig.~\ref{Fig2}.
The bottom and top lines correspond to all fullerene graphs and to IPR fullerene graphs,
respectively.
Explicit values of $C_n$ and quantity $g_n$ of the
graphs are presented in Table~\ref{TDistrAll}.
Since the minimal $g_n$ is equal to $9$,
we can formulate the following statement.

\begin{proposition}
There do not exist transmission irregular fullerene graphs with $n \le 216$ vertices.
\end{proposition}

Since the almost all fullerene graphs have no symmetries,
we believe that transmission irregular graphs exist for large
number of vertices.

\begin{problem}
Does there exist a transmission irregular fullerene graph (IPR fullerene graph)?
If yes, then what is the smallest order of such graphs?
\end{problem}

\clearpage

\begin{table}[t!h!]
\centering
\caption{Distribution of fullerene graphs of $F_{196}$ with respect to the Wiener \\ complexity $C_W$
         ($N$ is the number of graphs).} \label{TDistr100}
\footnotesize
\begin{tabular}{rr|rr|rr|rr|rr|rr} \hline
$C_W$&\ $N$&$C_W$ & $N$   &  $C_W$& $N$\ \ \ \  &  $C_W$ &$N$ \ \ \ \ & $C_W$ &$N$\ \ \ \   \\ \hline
 13 &    1 & 45 &  162    & 73 &  4635   & 101  & 289459   & 129 & 3966428  & 157& 1787105    \\
 14 &    2 & 46 &  176    & 74 &  5426   & 102  & 333904   & 130 & 4148282  & 158& 1534799    \\
 15 &    3 & 47 &  178    & 75 &  6367   & 103  & 380828   & 131 & 4323678  & 159& 1296521    \\
 16 &    1 & 48 &  134    & 76 &  7143   & 104  & 437958   & 132 & 4482070  & 160& 1080839   \\
 19 &    1 & 49 &  101    & 77 &  8442   & 105  & 498054   & 133 & 4630277  & 161& 885355     \\
 21 &    1 & 50 &   86    & 78 &  9680   & 106  & 564113   & 134 & 4766726  & 162& 713786     \\
 23 &    5 &  51&    78   & 79 & 10834   & 107  & 637114   & 135 & 4889247  & 163& 564259     \\
 24 &    3 &  52&    90   & 80 & 12451   & 108  & 718248   & 136 & 4999476  & 164& 438608     \\
 25 &    1 &  53&   114   & 81 & 13990   & 109  & 804261   & 137 & 5082927  & 165& 331636     \\
 26 &    4 &  54&   100   & 82 & 16160   & 110  & 899361   & 138 & 5147054  & 166& 246578     \\
 27 &    7 &  55&   112   & 83 & 18120   & 111  & 1000968  & 139 & 5186406  & 167& 178749    \\
 28 &    1 &  56&   132   & 84 & 20406   & 112  & 1110963  & 140 & 5195292  & 168& 126604    \\
 29 &    4 &  57&   173   &85  & 23226   & 113  & 1231994  & 141 & 5177020  & 169& 87271     \\
 30 &    4 &  58&   247   &86  & 26153   & 114  & 1357031  & 142 & 5131736  & 170& 58523     \\
 31 &    4 &  59&    268  &87  & 29857   & 115  & 1490391  & 143 & 5052595  & 171& 38065     \\
 32 &    8 &  60&    325  &88  & 34232   & 116  & 1632573  & 144 & 4943731  & 172& 23910     \\
 33 &   10 &  61&    429  &89  & 39558   & 117  & 1783764  & 145 & 4809637  & 173& 14592     \\
 34 &    8 &  62&    551  &90  & 46466   & 118  & 1941577  & 146 & 4643527  & 174&  8433     \\
 35 &   13 &  63&    619  & 91 & 55187   & 119  & 2108524  & 147 & 4450101  & 175&  4630     \\
 36 &   12 &  64&    846  & 92 & 65082   & 120  & 2278437  & 148 & 4236185  & 176&  2549     \\
 37 &   33 &  65&   1039  & 93 & 77669   & 121  & 2459465  & 149 & 4005961  & 177&  1318     \\
 38 &   27 &  66&   1268  & 94 & 92772   & 122  & 2636351  & 150 & 3750830  & 178&   653     \\
 39 &   48 &  67&   1587  & 95 & 110504  & 123  & 2825077  & 151 & 3479586  & 179&   306     \\
 40 &   60 &  68&   1777  & 96 & 130842  & 124  & 3016435  & 152 & 3198568  & 180&   130     \\
 41 &   78 &  69&   2267  & 97 & 155105  & 125  & 3210085  & 153 & 2912936  & 181&    71     \\
 42 &  121 &  70&   2704  & 98 & 181583  & 126  & 3401907  & 154 & 2624386  & 182&    26     \\
 43 &  132 &  71&   3279  & 99 & 214088  & 127  & 3594118  & 155 & 2339326  & 183&     8     \\
 44 &  153 &  72&   4016  &100 & 249142  & 128  & 3784693  & 156 & 2059994  & 184&     5     \\  \hline
\end{tabular}
\end{table}

\clearpage

\begin{table}[th]
\centering
\caption{Maximal Wiener complexity $C_n$ and $g_n=n-C_n$ of fullerene \\ graphs of $F_n$
         ($N$ is the number of graphs with $C_n$).} \label{TDistrAll}
\footnotesize
\begin{tabular}{cccc|cccc|cccc|cccc}
\hline
$n$ &$C_n$&$g_n$&$N$&   $n$&$C_n$&$g_n$&$N$&    $n$&$C_n$&$g_n$&$N$&  $n$&$C_n$&$g_n$&$N$ \\ \hline
 20 &  1 & 19  & 1 & 72  & 56  & 16 & 6 &  122  & 109 & 13 & 5 &172  & 160 & 12 & 1 \\
 24 &  2 & 22  & 1 & 74  & 61  & 13 & 1 &  124  & 111 & 13 & 5 &174  &164  & 10 & 1 \\
 26 &  2 & 24  & 1 & 76  & 63  & 13 & 1 &  126  & 115 & 11 & 1 &176  &165  & 11 & 1 \\
 28 &  5 & 23  & 1 & 78  & 64  & 14 & 2 &  128  & 117 & 11 & 1 &178  &167  & 11 & 2  \\
 30 &  7 & 23  & 1 & 80  & 66  & 14 & 2 &  130  & 118 & 12 & 3 &180  &168  & 12 & 1  \\
 32 &  9 & 23  & 1 & 82  &  71 &  11& 1 &  132  & 121 & 11 & 1 &182  &171  & 11 & 1  \\
 34 & 10 & 24  & 2 & 84  &  70 &  14& 2 &  134  & 123 & 11 & 3 &184  &172  & 12 & 4 \\
 36 & 14 & 22  & 1 & 86  &  73 &  13& 3 &  136  & 124 & 12 & 1 &186  &177  &  9 & 1  \\
 38 & 18 & 20  & 1 & 88  &  73 &  15& 7 &  138  & 127 & 11 & 2 &188  &177  & 11 & 1 \\
 40 & 19 & 21  & 1 & 90  & 79  & 11 & 1 &  140  & 131 &  9 & 1 &190  &180  & 10 & 1 \\
 42 & 22 & 20  & 1 & 92  & 80  & 12 & 1 &  142  & 131 & 11 & 1 &192  &181  & 11 & 1 \\
 44 & 25 & 19  & 1 & 94  & 82  & 12 & 1 &  144  & 132 & 12 & 2 &194  &183  & 11 & 2  \\
 46 & 25 & 21  & 4 & 96  & 84  & 12 & 2 &  146  & 134 & 12 & 4 &196  &184  & 12 & 5 \\
 48 & 30 & 18  & 1 & 98  & 86  & 12 & 1 &  148  & 136 & 12 & 1 &198  &187  & 11 & 2 \\
 50 & 35 & 15  & 1 & 100 & 89  & 11 & 1 &  150  & 138 & 12 & 4 &200  &189  & 11 & 2 \\
 52 & 36 & 16  & 1 & 102 & 89  & 13 & 4 &  152  & 141 & 11 & 1 &202  &193  &  9 & 1  \\
 54 & 37 & 17  & 1 & 104 & 91  & 13 & 3 &  154  & 144 & 10 & 1 &204  &192  & 12 & 1 \\
 56 & 40 & 16  & 1 & 106 & 93  & 13 & 3 &  156  & 144 & 12 & 3 &206  &195  & 11 & 3 \\
 58 & 43 & 15  & 2 & 108 & 96  & 12 & 1 &  158  & 147 & 11 & 1 &208  &198  & 10 & 1 \\
 60 & 44 & 16  & 3 & 110 & 98  & 12 & 2 &  160  & 148 & 12 & 2 &210  &199  & 11 & 1 \\
 62 & 46 & 16  & 3 & 112 & 100 & 12 & 1 &  162  & 151 & 11 & 2 &212  &199  & 13 & 8 \\
 64 & 49 & 15  & 5 & 114 & 102 & 12 & 1 &  164  & 153 & 11 & 1 &214  &202  & 12 & 3 \\
 66 & 50 & 16  & 2 & 116 & 102 & 14 & 4 &  166  & 155 & 11 & 2 &216  &204  & 12 & 4 \\
 68 & 56 & 12  & 1 & 118 & 106 & 12 & 3 &  168  & 157 & 11 & 2 & -   &     &    &   \\
 70 & 56 & 14  & 1 & 120 & 110 & 10 & 2 &  170  & 159 & 11 & 2 & -   &     &    &   \\ \hline \end{tabular}
\end{table}

\section{Graphs with the maximal Wiener complexity}

In this section, we study the following problem:
can the Wiener index of a fullerene graph with the maximal Wiener complexity be maximal?
Numerical data for the Wiener indices of fullerene graphs of order $n \le 216$
are presented in Table~\ref{TWiener}.
Here three columns $C_n$, $W$, and $D$ are the maximal Wiener complexity,
the Wiener index and the diameter of graphs with $C_n$, respectively.
Three columns $W_m$, $C_W$, and $D$ contain the maximal Wiener index,
the Wiener complexity and the diameter of graphs with $W_m$.

Based on data of Tables~\ref{TDistrAll} and \ref{TWiener}, one can make the following observations.
\begin{itemize}
\item
Several fullerene graphs of fixed $n$ may have the maximal Wiener complexity $C_n$
while the only one fullerene graph has the maximal Wiener index.

\item
Wiener indices of fullerene graphs with fixed $C_n$ ($|\, F_n| > 1$) are not maximal except
graphs of order $n=28$ with $W=1198$ ($|\, F_{28}| = 2$).

\item
Almost all fullerene graphs with fixed $C_n$ have distinct Wiener indices.
The only exception are graphs of order 46 with $W=4289$
(the sequences of their vertex transmissions are distinct).
\end{itemize}

\pagebreak

\begin{table}[t!h!]
\centering
\caption{\normalsize Maximal Wiener complexity and Wiener indices of fullerene graphs.} \label{TWiener}
\footnotesize
\begin{tabular}{|r@{\hspace{2mm}}|r@{\hspace{2mm}}r@{\hspace{3mm}}r|r@{\hspace{2mm}}r@{\hspace{2mm}}r@{\hspace{2mm}}r
                |@{\hspace{1mm}}r@{\hspace{1mm}}|
                @{\hspace{1mm}}r@{\hspace{2mm}}|rr@{\hspace{3mm}}r|r@{\hspace{2mm}}r@{\hspace{2mm}}r@{\hspace{2mm}}r|}
\hline
$n$ & $C_n$&$W$    &$D$  & $W_m$ &$C_W$ & $D$&$t$& &$n$ &$C_n$&$W$    &$D$   &$W_m$ &$C_W$ &$D$& $t$ \\ \hline
 20 &   1 & 500   & 5   & 500   & 1  & 5  &    & & 84 & 70 & 19939 & 13   & 21754 & 21& 15&$c1$ \\
 24 &   2 & 804   & 5   & 804   & 2  & 5  &    & &    &    & 20076 & 13   &       &   &   &   \\
 26 &   2 & 987   & 6   & 987   & 2  & 6  &$b$ & & 86 & 73 & 21404 & 13   & 23467& 8  & 16&$b$ \\
 28 &   5 & 1198  & 6   & 1198  & 5  & 6  &    & &    &    & 21521 & 13   &      &    &   &     \\
 30 &   7 & 1431  & 6   & 1435  & 3  & 6  &$a$ & &    &    & 21593 & 13   &      &    &   &     \\
 32 &   9 & 1688  & 6   & 1696  & 3  & 7  &$b$ & & 88 & 73 & 22359 & 13   &24714 & 21 & 16& $c2$ \\
 34 &  10 & 1973  & 7   & 1978  & 10 & 7  &    & &    &    & 22421 & 13   &      &    &   &   \\
    &     & 1978  & 7   &       &    &    &    & &    &    & 22604 & 13   &      &    &   &  \\
 36 &  14 & 2288  & 7   & 2298  & 8  & 7  &$c1$ & &    &    & 22616 & 13   &      &    &   &  \\
 38 &  18 & 2627  & 7   & 2651  & 4  & 8  &$b$ & &    &    & 22619 & 13   &      &    &   &  \\
 40 &  19 & 3001  & 7   & 3035  & 4  & 8  &$a$ & &    &    & 22750 & 14   &      &    &   &   \\
 42 &  22 & 3397  & 8   & 3415  & 19 & 8  &$d1$ & &    &    & 22939 & 14   &      &    &   &  \\
 44 &  25 & 3830  & 8   & 3888  & 4  & 9  &$b$ & &90  & 79 & 23923 & 14   &27155 &  9 & 17&$a$ \\
 46 &  25 & 4285  & 8   & 4322  & 19 & 9  &$d2$ & &92  & 80 & 25731 & 15   &28256 &  8 & 17&$b$ \\
    &     & 4289  & 8   &       &    &    &    & &94  & 82 & 26793 & 14   &28910 & 44 & 17&$d2$ \\
    &     & 4289  & 8   &       &    &    &    & &96  & 84 & 28274 & 14   &31418 & 24 & 17&$c1$ \\
    &     & 4291  & 8   &       &    &    &    & &    &    & 28317 & 15   &      &    &   &  \\
 48 &  30 & 4795  & 9   & 4858  & 12 & 9  &$c1$ & &98  & 86 & 30068 & 15   &33651 &  9 & 18&$b$ \\
 50 &  35 & 5310  & 9   & 5455  & 5  & 9  &$a$ & &100 &89  & 31196 & 15   &36580 & 10 & 19&$a$ \\
 52 &  36 & 5876  & 9   & 5994  & 13 & 10 &$c2$ & &102 &89  & 32984 & 15   &36206 & 47 & 18&$d1$ \\
 54 &  37 & 6475  & 9   & 6558  & 22 & 10 &$d1$ & &    &    & 33070 & 15   &      &    &   &  \\
 56 &  40 & 7114  & 10  & 7352  & 5  & 11 &$b$ & &    &    & 33226 & 15   &      &    &   &  \\
 58 &  43 & 7782  & 10  & 7910  & 25 & 11 &$d2$ & &    &    & 33505 & 15   &      &    &   &  \\
    &     & 7822  & 10  &       &    &    &    & &104 &91  & 34402 & 15   &39688 & 9  & 19&$b$ \\
 60 &  44 & 8437  & 10  & 8880  & 6  & 11 &$a$ & &    &    & 34529 & 15   &      &    &   &  \\
    &     & 8466  & 10  &       &    &    &    & &    &    & 36801 & 17   &      &    &   &  \\
    &     & 8490  & 10  &       &    &    &    & &106 &93  & 36648 & 16   &40278 & 47 & 19&$d2$ \\
 62 &  46 & 9202  & 10  & 9651  & 6  & 12 &$b$ & &    &    & 36664 & 16   &      &    &   &  \\
    &     & 9220  & 11  &       &    &    &    & &    &    & 37594 & 17   &      &    &   &  \\
    &     & 9250  & 11  &       &    &    &    & &108 &96  & 38033 & 15   & 43578& 27 & 19&$c1$ \\
 64 &  49 & 9988  & 11  & 10410 & 15 & 12 &$c2$ & &110 &98  & 40154 & 16   & 48005& 11 & 21&$a$ \\
    &     & 9993  & 11  &       &    &    &    & &    &    & 41419 & 17   &      &    &   &   \\
    &     &10003  & 11  &       &    &    &    & &112 &100 & 41940 & 17   &48234 & 27 & 20&$c2$ \\
    &     &10013  & 11  &       &    &    &    & &114 &102 & 43885 & 16   &49318 & 52 & 20&$d1$ \\
    &     &10016  & 11  &       &    &    &    & &116 &102 & 45437 & 16   &53832 & 10 & 21&$b$ \\
 66 &  50 &10814  & 11  &11126  & 30 & 12 &$d1$ & &    &    & 46632 & 17   &      &    &   &  \\
    &     &10842  & 11  &       &    &    &    & &    &    & 46798 & 17   &      &    &   &  \\
 68 &  56 &11714  & 11  &12376  & 6  & 13 &$b$ & &    &    & 47927 & 18   &      &    &   &  \\
 70 &  56 &12589  & 11  &13505  & 7  & 13 &$a$ & &118 &106 & 47059 & 15   &54310 & 50 & 21&$d2$ \\
 72 &  56 &13407  & 11  &14298  & 18 & 13 &$c1$ & &    &    & 47489 & 16   &      &    &   &  \\
    &     & 13448 & 11  &       &    &    &    & &    &    & 47697 & 16   &      &    &   &  \\
    &     & 13453 & 11  &       &    &    &    & &120 &110 & 49143 & 16   & 61630& 12 & 23&$a$ \\
    &     & 13567 & 12  &       &    &    &    & &    &    & 49606 & 17   &      &    &   &  \\
    &     & 13578 & 12  &       &    &    &    & &122 &109 & 51344 & 16   & 62011& 11 & 22&$b$ \\
    &     & 13766 & 12  &       &    &    &    & &    &    & 51456 & 16   &       &   &   &  \\
 74 &  61 & 14521 & 12  & 15563 & 7  & 14 &$b$ & &    &    & 51933 & 17   &       &   &   &  \\
 76 &  63 & 15867 & 13  & 16554 & 18 & 14 &$c2$ & &    &    & 52974 & 17   &       &   &   &  \\
 78 &  64 & 16834 & 13  & 17398 & 37 & 14 &$d1$ & &    &    & 53070 & 17   &       &   &   &  \\
    &     & 16877 & 13  &       &    &    &    & &124 &111 & 54105 & 17   & 64170 & 30& 22&$c2$ \\
 80 &  66 & 17727 & 13  & 19530 & 8  & 15 &$a$ & &    &    & 55050 & 18   &       &   &   &  \\
    &     & 17832 & 13  &       &    &    &    & &    &    & 55789 & 18   &       &   &   &  \\
 82 &  71 & 19075 & 13  & 19918 & 38 & 15 &$d2$ & &    &    & 57358 & 19   &       &   &   &  \\
    &     &       &     &       &    &    &    & &    &    & 57473 & 19   &       &   &   &  \\ \hline
\end{tabular}
\end{table}

\clearpage

\setcounter{table}{2}

\begin{table}[ht]
\centering
\caption{\normalsize Maximal Wiener complexity and Wiener indices of fullerene graphs (\emph{continue}).}
\footnotesize
\begin{tabular}{|r@{\hspace{2mm}}|r@{\hspace{2mm}}r@{\hspace{2mm}}r|r@{\hspace{2mm}}r@{\hspace{2mm}}r@{\hspace{2mm}}r
               |@{\hspace{1mm}}r@{\hspace{1mm}}|
                @{\hspace{1mm}}r@{\hspace{2mm}}|rr@{\hspace{2mm}}r|r@{\hspace{2mm}}r@{\hspace{2mm}}r@{\hspace{2mm}}r|}
\hline
$n$ &  $C_m$&$W$    &$D$  & $W_m$ &$C$ & $D$&$t$& &$n$ &$C_m$&$W$    &$D$   &$W_m$ &$C$ &$D$& $t$ \\ \hline
 126 & 115  & 57238 & 18  & 65286 & 57 & 22 &$d1$& &178 &  167 & 141743& 23   &174510& 65 &31 &$d2$  \\
 128 & 117  & 60434 & 18  & 70976 & 11 & 23 &$b$& &    &      & 150696& 25   &      &    &   &     \\
 130 &  118 & 63736 & 19  & 77655 & 13 & 25 &$a$& &180 &  168 & 139697& 21   &200780& 18 &35 &$a$  \\
     &      & 63922 & 19  &       &    &    &   & &182 &  171 & 144410& 21   &192971& 16 & 32&$b$  \\
     &      & 65396 & 20  &       &    &    &   & &184 &  172 & 146581& 22   &197130& 45 & 32&$c2$  \\
 132 &  121 & 62917 & 17  & 76538 & 33 & 23 &$c1$& &    &      & 147054& 21   &      &    &   &     \\
 134 &  123 & 64935 & 17  & 80763 & 12 & 24 &$b$& &    &      & 153615& 23   &      &    &   &    \\
     &      & 65225 & 17  &       &    &    &   & &    &      & 154923& 23   &      &    &   &     \\
     &      & 68161 & 19  &       &    &    &   & &186 &  177 & 167300& 26   &198046& 82 & 32&$d1$  \\
 136 &  124 & 69838 & 19  & 83274 & 33 & 24 &$c2$& &188 &  177 & 154868& 21   &211776& 16 & 33&$b$  \\
 138 &  127 & 72311 & 19  & 84398 & 62 & 24 &$d1$& &190 &  180 & 169849& 24   &235405& 19 & 37&$a$  \\
     &      & 73771 & 19  &       &    &    &   & &192 &  181 & 163370& 22   &222778& 48 & 33&$c1$  \\
 140 &  131 & 73644 & 19  & 96280 & 14 & 27 &$a$& &194 &  183 & 187947& 27   &231763& 17 & 34&$b$  \\
 142 &  131 & 79852 & 20  & 91518 & 56 & 25 &$d2$& &    &      & 191290& 27   &      &    &   &     \\
 144 &  132 & 77934 & 18  & 97914 & 36 & 25 &$c1$& &196 &  184 & 174774& 23   &236394& 48 & 34&$c2$  \\
     &      & 78924 & 18  &       &    &    &   & &    &      & 178192& 24   &      &    &   &     \\
 146 &  134 & 86095 & 20  &102947 & 13 & 26 &$b$& &    &      & 178529& 25   &      &    &   &     \\
     &      & 86287 & 21  &       &    &    &   & &    &      & 179284& 25   &      &    &   &    \\
     &      & 87298 & 21  &       &    &    &   & &    &      & 184011& 24   &      &    &   &    \\
     &      & 87442 & 21  &       &    &    &   & &198 &  187 & 177296& 23   &237198& 87 & 34& $d1$ \\
 148 &  136 & 86432 & 20  &105834 & 36 & 26 &$c2$& &    &      & 189530& 25   &      &    &   &     \\
 150 &  138 & 87886 & 19  &117705 & 15 & 29 &$a$& &200 &  189 & 180683& 23   &273830& 20 & 39& $a$ \\
     &      & 88860 & 19  &       &    &    &   & &    &      & 192365& 25   &      &    &   &     \\
     &      & 92732 & 21  &       &    &    &   & &202 &  193 & 210388& 28   &251318& 71 & 35& $d2$  \\
     &      & 93898 & 21  &       &    &    &   & &204 &  192 & 189450& 22   &265274& 51 & 35& $c1$  \\
 152 &  141 & 92988 & 20  &115416 & 13 & 27 &$b$& &206 &  195 & 221909& 28   &275427& 18 & 36&  $b$ \\
 154 &  144 & 97359 & 21  &115270 & 59 & 27 &$d2$& &    &      & 221995& 28   &      &    &   &     \\
 156 &  144 & 95579 & 19  &122938 & 39 & 27 &$c1$& &    &      & 222097& 28   &      &    &   &     \\
     &      & 96997 & 19  &       &    &    &   & &208 &  198 &201644 & 23   &280554& 51 & 36& $c2$  \\
     &      & 98864 & 21  &       &    &    &   & &210 &  199 &238572 & 29   &316255& 21 & 41& $a$  \\
 158 &  147 &100055 & 19  &128851 & 14 & 28 &$b$& &212 &  199 &207617 & 23   &299176& 18 & 37& $b$ \\
 160 &  148 &103952 & 20  &142130 & 16 & 31 &$a$& &    &      &207975 & 23   &      &    &   &    \\
     &      &108170 & 22  &      &     &    &   & &    &      &209707 & 23   &      &    &   &   \\
 162 &  151 & 104909& 19   &133206& 72 & 28&$d1$ & &    &      &211942 & 24   &      &    &   &    \\
 &      & 116278& 23   &      &    &   &    & &    &      &215779 & 24   &      &    &   &    \\
 164 &  153 & 110088& 20   &143288& 14 & 29&$b$ & &    &      &228285 & 26   &      &    &   &     \\
 166 &  155 & 117531& 21   &142838& 62 & 29&$d2$ & &    &      &228507 & 26   &      &    &   &     \\
     &      & 119485& 23   &      &    &   &    & &    &      &228922 & 26   &      &    &   &    \\
 168 &  157 & 114316& 18   &151898& 42 & 29&$c1$ & &214 &  202 &226652 & 25   &297030& 74 & 37& $d2$ \\
     &      & 126632& 23   &      &    &   &    & &    &      &247352 & 29   &      &    &   &     \\
 170 &  159 & 123193& 22   &169755& 17 & 33&$a$ & &    &      &250978 & 29   &      &    &   &     \\
     &      & 130548& 24   &      &    &   &    & &216 &  204 &220131 & 23   &312858& 54 & 37& $c1$ \\
 172 &  160 & 129708& 22   &162474& 42 & 30&$c2$ & &    &      &226928 & 25   &      &    &   &     \\
 174 &  164 & 131354& 23   &163478& 77 & 30&$d1$ & &    &      &240920 & 27   &      &    &   &     \\
 176 &  165 &130105 & 20   &175312& 15 &31 &$b$ & &    &      &270770 & 33   &      &    &   & \\  \hline
\end{tabular}
\end{table}

\clearpage

\begin{itemize}
\item
The diameter of graphs with fixed $C_n$ are not maximal for $n \ge 52$.

\item
Fullerene graphs with the maximal Wiener index have the maximal diameter.
The values of the Wiener complexity $C_W$ can vary greatly.
This can be partially explained by the appearance of symmetries in graphs.
\end{itemize}

It is of interest how the pentagons are distributed among hexagons
for fullerene graphs with the maximal Wiener complexity
(see Tables~\ref{TDistrAll} and \ref{TWiener}).
Does there exist any regularity in the distribution of pentagons?
Table~\ref{TConn} gives some
information on the occurrence of pentagonal parts of a particular size.
Here $N$ is the number of graphs in which pentagons form $N_p$ isolated connected parts.

\vspace{1mm}

\begin{table}[h]
\centering
\caption{The number of graphs with $N_p$ isolated pentagonal parts.}  \label{TConn}
\begin{tabular}{l|rrrrrrrr} \hline
 $N_p$  &  1 & 2 & 3  & 4  & 5  & 6  & 7 & 8  \\
 $N$    & 9  & 8 & 27 & 61 & 42 & 40 & 7 & 1 \\ \hline
\end{tabular}
\end{table}

Table~\ref{TIsolC5} shows how many fullerene graphs with the maximal Wiener complexity
have isolated pentagons (an isolated pentagon forms a part).
Here $N$ is the number of graphs having  $N_5$ isolated pentagons.
Does there exist an IPR fullerene graph with maximal Wiener complexity $C_n$
(lines of Fig.~\ref{Fig2} will have intersection)?

\vspace{1mm}

\begin{table}[th]
\centering
\caption{The number graphs of with $N_5$ isolated pentagons.}  \label{TIsolC5}
\begin{tabular}{l|rrrrrr}  \hline
 $N_5$  & 0 &  1  & 2 & 3  & 4  & 5   \\
 $N$    & 23&  56 & 44& 44 & 24 & 4  \\ \hline
\end{tabular}
\end{table}

\section{Graphs with the maximal Wiener index}

\begin{figure}[th]
\vspace*{5mm}
\center{\includegraphics[width=0.8\linewidth]{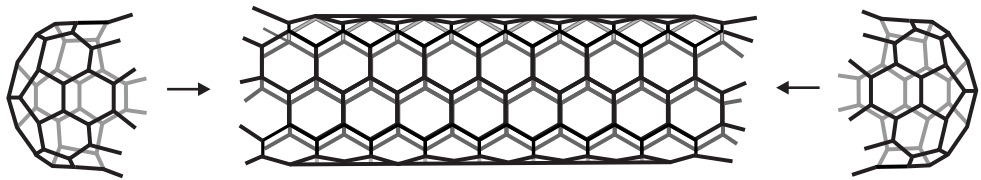}}
\caption{Construction of a nanotubical fullerene graph with two caps.}
\vspace*{2mm}
\label{Fig3}
\end{figure}

Wiener index of fullerene graphs are studied in
\cite{Aliz14,Ando17,Ashr11,Fowl01,Ghor13-1,Ghor17,Ghor13-2,Ghos18,Grao11,Hua14,Iran09}.
There is a class of fullerene graphs of tubular shapes,
called \emph{nanotubical fullerene graphs}.
They are cylindrical in shape, with the two ends capped
by  a subgraph containing six pentagons and possibly some hexagons
called \emph{caps} (see an illustration in Fig.~\ref{Fig3}).

\begin{figure}[t]
\center{\includegraphics[width=0.8\linewidth]{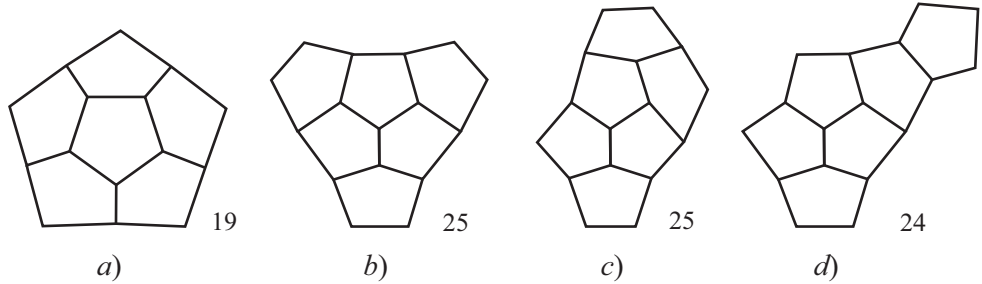}}
\caption{Caps for nanotubical fullerene graphs with the maximal Wiener index.}
\label{Fig4}
\end{figure}

Consider fullerene graphs with the maximal Wiener indices
(see Table~\ref{TWiener}).
Five graphs of $F_{20}$--$F_{28}$ and $F_{34}$ contain one pentagonal part
and other 93 graphs possess two pentagonal parts.
Two pentagonal parts of every fullerene graph
are the same and contain diametrical vertices.
Therefore such graphs are nanotubical fullerene
graphs with caps containing identical pentagonal parts.
All types of such parts are depicted in Fig.~\ref{Fig4}.
The number of fullerene graphs having a given part is shown near diagrams.
A type of a cap is determined by the type of its pentagonal part.
Types of caps of fullerene graphs are presented
by the corresponding notation in column $t$ of Table~\ref{TWiener}.
Constructive approaches for enumeration of various caps
were proposed in
\cite{Brin02,Brin99}.
Consider every kind of cap types.

\begin{figure}[th]
\center{\includegraphics[width=0.8\linewidth]{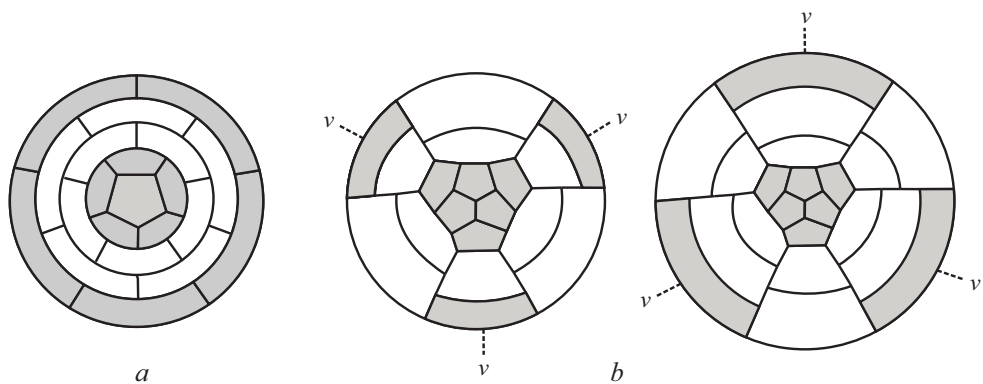}}
\caption{Structure of fullerene graphs with caps of types $a$ and $b$.}
\label{Fig5}
\end{figure}

1. \emph{Type $a$}. Caps of type $a$ define so-called (5,0)-nanotubical
fullerene graphs.
The structure of graphs of this infinite family $T_a$ is clear
from an example in Fig.~\ref{Fig5}a.
Diameter and the Wiener index of such fullerene graphs were studied in
\cite{Aliz14}.
To indicate the order of graph $G$, we will use notation $G_n$.

\vspace{-2mm}

\begin{proposition} {\rm \cite{Aliz14}}  \label{Wa}
Let $G_n$ be a nanotubical fullerene graph with caps of type $a$.
It has $n=10k$ vertices, $k \ge 2$. Then $C_W(G_n)=k$, $D(G_n)=2k-1$,  and
$W(G_{20}) = 500$, $W(G_{30}) = 1435$, $W(G_{40}) = 3035$, and for $n \ge 50$,
\begin{eqnarray*}
W(G_n) & = & \frac{1}{30} \left( n^3 + 1175 n - 20100 \right).
\end{eqnarray*}
\end{proposition}

Based on numerical data of Table~\ref{TWiener}, the similar results have been
obtained for fullerene graphs
of order $n \le 216$ with caps of the other three types.

\newpage

2. \emph{Type $b$}. The structure of graphs of the corresponding family $T_b$ is clear from examples
of Fig.~\ref{Fig5}b. Vertices marked by $v$ should be identified
in every graph. Table~\ref{TWiener} contains 26 such graphs.

\begin{proposition} \label{Wb}
Let $G_n$ be a nanotubical fullerene graph with caps of type $b$.
It has $n=6k-4$ vertices, $k \ge 5$. Then
$C_W(G_n)=\lceil k/2 \rceil$, $D(G_n)=k+1$, and for $n \ge 26$,
\begin{eqnarray*}
W(G_n) & = & \frac{1}{36} \left( n^3 + 27 n^2 + 156 n - 4352  \right).
\end{eqnarray*}
\end{proposition}

Two caps of  type $b$ have adjacent pentagonal rings only for $k=5$ .
If fullerenes with caps of types $a$ and $b$ have the same number
of faces ($n=10k$), then the graph with caps
of type $a$ has the maximal Wiener index.

\vspace{1mm}

3. \emph{Type $c$}. Fullerene graphs with caps of type $c$ will be splitted
into disjoint fa\-mi\-lies,
$T_c = T_{c1} \cup T_{c2}$. The corresponding graphs are
marked in column $t$ of Table~\ref{TWiener}
by $c1$ (13 graphs) and $c2$ (12 graphs).
The numbers of vertices of graphs are
given in Table~\ref{Tcd}.
The orders of graphs of $T_c$  do not coincide with the orders of
graphs from the set $T_a \cup T_b$.

\begin{proposition} \label{Wc}
a) Let $G_n$ be a nanotubical fullerene graph of family $T_{c1}$. Then for $n \ge 36$,
\begin{eqnarray*}
W(G_n) & = & \frac{1}{36} \left( n^3 + 24 n^2 + 336 n - 7128  \right).
\end{eqnarray*}
The Wiener complexity and the diameter of $G_n$ are shown in Table~\ref{Tcd}.
One value should be corrected for $k=0$ (see a cell of Table~\ref{Tcd} with
mark *): $C_W(G_{36})=8$ instead of $9$.

b) Let $G_n$ be a nanotubical fullerene graph of family $T_{c2}$. Then for $n \ge 52$,
\begin{eqnarray*}
W(G_n) & = & \frac{1}{36} \left( n^3 + 24 n^2 + 336 n - 7192  \right).
\end{eqnarray*}
The Wiener complexity and the diameter of $G_n$ are shown in Table~\ref{Tcd}.
One value should be corrected for $k=0$: $C_W(G_{52})=13$ instead of $12$.
\end{proposition}

\begin{table}[t]
\centering
\caption{Parameters of fullerene graphs with $n \le 216$ vertices
 and caps of types $c$ and $d$. Here $k\ge 0$ for all expressions.} \label{Tcd}
%\footnotesize
\begin{tabular}{lll|lll} \hline
\multicolumn{3}{c|}{family $T_{c1}$}  & \multicolumn{3}{c}{family $T_{c2}$} \\ \hline
\ \ \ \ \  $n$ &\ \ \ $C_W$&\ \ \ \ $D$&
\ \ \ \ \  $n$ &\ \ \ $C_W$&\ \ \ \ $D$  \\  \hline
$60k+36$ & $15k+9^{\,*}$ & $10k+7$  & $60k+76$ & $15k+18$       & $10k+14$  \\
$60k+48$ & $15k+12$      & $10k+9$  & $60k+88$ & $15k+21$       & $10k+16$  \\
$60k+72$ & $15k+18$      & $10k+13$ & $60k+52$ & $15k+12^{\,*}$ & $10k+10$  \\
$60k+84$ & $15k+21$      & $10k+15$ & $60k+64$ & $15k+15$       & $10k+12$  \\ \hline
\multicolumn{3}{c|}{family $T_{d1}$} & \multicolumn{3}{c}{family $T_{d2}$} \\ \hline
\ \ \ \ \  $n$ &\ \ \ $C_W$&\ \ \ \ $D$&
\ \ \ \ \  $n$ &\ \ \ $C_W$&\ \ \ \ $D$  \\  \hline
$60k+66$ & $25k+32^{\,*}$ & $10k+2$  & $60k+106$& $15k+47$       & $10k+19$  \\
$60k+78$ & $25k+37$       & $10k+4$  & $60k+58$ & $15k+35^{\,*}$ & $10k+11$  \\
$60k+102$& $25k+47$       & $10k+8$  & $60k+82$ & $15k+41^{\,*}$ & $10k+15$  \\
$60k+54$ & $25k+27^{\,*}$ & $10k$    & $60k+94$ & $15k+44$       & $10k+17$  \\ \hline
\end{tabular}
\end{table}

\vspace{1mm}

4. \emph{Type $d$}. Fullerene graphs with caps of type $d$ will be also splitted into
two disjoint families, $T_d = T_{d1} \cup T_{d2}$.
 The both families have 12 members
 (see graphs with marks $d1$ and $d2$ in column $t$ of Table~\ref{TWiener}).
The numbers of vertices of graphs of  $T_d$
are shown in Table~\ref{Tcd}.
The orders of graphs of $T_d$  do not coincide with the orders of
graphs from the set $T_a \cup T_b \cup T_c$.

\begin{proposition} \label{Wd}
a) Let $G_n$ be a nanotubical fullerene graph of family $T_{d1}$. Then
$W(G_{42})=3415$ and for $n \ge 54$,
\begin{eqnarray*}
W(G_n) & = & \frac{1}{36} \left( n^3 + 15 n^2 + 1068 n - 22788  \right).
\end{eqnarray*}
The Wiener complexity and the diameter of $G_n$ are shown in Table~\ref{Tcd}.
Two values should be corrected for $k=0$ (see cells of Table~\ref{Tcd} with mark *):
$C_W(G_{66})=30$  instead of $32$ and $C_W(G_{54})=22$ instead of $27$.

b) Let $G_n$ be a nanotubical fullerene graph of family $T_{d2}$.
Then $W(G_{46})=4322$ and for $n \ge 58$,
\begin{eqnarray*}
W(G_n) & = & \frac{1}{36} \left( n^3 + 15 n^2 + 1068 n - 22756  \right).
\end{eqnarray*}
The Wiener complexity and the diameter of $G_n$ are shown in Table~\ref{Tcd}.
Two values should be corrected for $k=0$:
$C_W(G_{58})=25$  instead of $35$ and $C_W(G_{82})=38$ instead of $41$.
\end{proposition}

The above considerations of fullerene graphs with $n \le 216$ vertices lead to the following conjectures
for all fullerene graphs.

\begin{conjecture}
If a fullerene graph of an arbitrary order has the maximal Wiener index, then it is a nanotubical fullerene graph with caps
of types $a$--$d$ and its Wiener index is given by Propositions \ref{Wa}--\ref{Wd}.
\end{conjecture}

\begin{conjecture}
The Wiener complexity and the diameter of fullerene graphs of an arbitrary order
having the maximal Wiener index are given in Propositions \ref{Wa}--\ref{Wd}.
\end{conjecture}

\section*{Acknowledgements}
This work was supported by the Laboratory of Topology and Dynamics,
Novosibirsk State University (contract no.~14.Y26.31.0025 with the Ministry
of Education and Science of the Russian Federation).

\end{document}